%% file: transitivity.tex
\documentclass[11pt,reqno]{amsart}
\usepackage{amsbsy,amsfonts,amsmath,amssymb,amscd,amsthm}
\usepackage{mathrsfs,float}
\usepackage[mathcal]{euscript}
\usepackage{subfigure,enumitem,graphicx,epstopdf}
\usepackage[a4paper,text={6.1in,8in},centering]{geometry}
\usepackage{pxfonts,epstopdf}
\usepackage[colorlinks=true,linkcolor=blue,citecolor=green]{hyperref}
\usepackage{srcltx}
\usepackage[margin=20pt,font=small,labelfont=bf,textfont=it]{caption}

\hypersetup{linktocpage}

\small\normalsize
\theoremstyle{plain}
\newtheorem{mainthm}{Theorem}

\newtheorem{thm}{Theorem}[section]
\newtheorem{cor}[thm]{Corollary}
\newtheorem{lem}[thm]{Lemma}

\newtheorem{prop}[thm]{Proposition}

\newtheorem{defi}[thm]{Definition}
\newtheorem{scho}[thm]{Scholium}
\theoremstyle{definition}

\newtheorem{rem}[thm]{Remark}

\newcommand{\eqdef}{\stackrel{\scriptscriptstyle\rm def}{=}}

\makeatletter
\def\l@part{\@tocline{0}{-2pt}{1pc}{}{}}
\def\l@section{\@tocline{1}{-2pt}{1pc}{4.6em}{}}
\renewcommand{\tocpart}[3]{%
  \indentlabel{\@ifnotempty{#2}{\makebox[2.3em][l]{%
    \ignorespaces#1 #2.\hfill}}}\bf{#3}}
\renewcommand{\tocsection}[3]{%
  \indentlabel{\@ifnotempty{#2}{\hspace*{2.3em}\makebox[2.3em][l]{%
    \ignorespaces#1 #2.\hfill}}}#3}

\makeatother

\begin{document}
\title[Robust non-hyperbolic transitivity]{Robustly non-hyperbolic transitive symplectic dynamics}
\author[Barrientos]{Pablo G. Barrientos}
\address{\centerline{Instituto de Matem\'atica e Estat\'istica, UFF}
    \centerline{Rua M\'ario Santos Braga s/n - Campus Valonguinhos, Niter\'oi,  Brazil}}
\email{pgbarrientos@id.uff.br}
\author[Raibekas]{Artem Raibekas}
\address{\centerline{Instituto de Matem\'atica e Estat\'istica, UFF}
   \centerline{Rua M\'ario Santos Braga s/n - Campus Valonguinhos, Niter\'oi,  Brazil}}
\email{artem@mat.uff.br}


\begin{abstract}
We construct symplectomorphisms in dimension $d\geq 4$ having a
semi-local robustly transitive partially hyperbolic set containing
$C^2$-robust homoclinic tangencies of any codimension $c$ with
$0<c\leq d/2$.
\end{abstract}

 \maketitle
\thispagestyle{empty}


\input  transitivity-intro

\input  transitivity-preliminaries
\input  transitivity-criteria
\input  transitivity-symplectic
\subsection*{Acknowledgements}
We are grateful to Carlos Menino Coton, Meysam Nassiri, and Michele Triestino for discussions and helpful suggestions. During the preparation of this
article the first author was supported by MTM2014-56953-P, and the second author by FAPERJ-INST 111.595/2014.


\bibliographystyle{abbrv}
\bibliography{br-bib}

\end{document}

%% file: transitivity-intro.tex
\section{Introduction}
According to  KAM theory, for symplectomorphisms close to an
integrable one, an orbit, with large probability, belongs to an
invariant torus and thus stays bounded for all time. Furthermore,
orbits which are close, take a large time to escape from a neighborhood of
these tori. However, in higher dimensional systems the action
variables may, a priori, exhibit considerable change showing
unbounded orbits~\cite{Ar64}. Such behavior was coined
with the term of \emph{Arnold diffusion} or \emph{instability}.
Diffusion orbits can be constructed by using normally hyperbolic
invariant laminations~\cite{Ll02}. But, this kind of invariant
laminations also allow the construction of robustly transitive sets which
give a more complex mechanism of diffusion~\cite{NP12}. In this
paper we continue the work of~\cite{NP12} and the study of these large (semi-local)
transitive sets and provide new examples in this context.

It was unknown if the examples constructed in~\cite{NP12} were
robustly non-hyperbolic. One of the classical tools to
create robustly non-hyperbolic dynamics is via a heterodimensional cycle~\cite{BD96}.
However in the symplectic case
this idea fails since all hyperbolic periodic points have
the same stability index and thus, there are no heterodimensional
cycles. The other classical approach to
destroy hyperbolicity is the construction of robust homoclinic tangencies~\cite{New70}. In the symplectic setting
known results on persistence of tangencies are restricted to
area-preserving diffeomorphisms~\cite{Du99,Du00,MR97}. Motivated
by this problem, we have developed in~\cite{BR17} a new
method to construct robust homoclinic tangencies in higher
dimensions. These techniques can be applied in the symplectic
framework as we show in this paper. As a consequence, we extend
the results of~\cite{NP12} showing that the semi-local transitive
sets can be made robustly non-hyperbolic.

%
%

Assume that $N$ and $M$ are symplectic connected manifolds (not
necessarily compact).  Let $F:N\to N$ be a $C^r$-symplectomorphism
having a hyperbolic set $\Lambda \subset N$ conjugated to a full
shift with a big enough set of symbols (depending only on
the dimension of $M$).  In order to state our main theorem we also
need the following notion:

A diffeomorphism $f$ has a \emph{homoclinic tangency of
codimension $c>0$} if there is a pair of points $P$ and $Q$
belonging to the same transitive hyperbolic set so that the
unstable invariant manifold of $P$ and the stable invariant
manifold of $Q$ have a non-transverse intersection $Y$ of
codimension $c$. That is,
$$
  Y \in W^u(P)\cap W^s(Q) \quad \text{and} \quad  c= \dim T_Y W^u(P)\cap T_Y W^s(Q).
$$
Now we are ready to state the main result:

\begin{mainthm}
\label{thmD} There is an arc $\{f_\epsilon\}_{\epsilon\geq 0}$ of
$C^r$-symplectomorphisms of $N \times M$  such that
 $f_0 = F \times \mathrm{id}$ and for $\varepsilon>0$,
 any small enough  $C^2$-perturbation
 $g$ of $f_\varepsilon$ has
 \begin{enumerate}[leftmargin=0.6cm, itemsep=0.1cm]
 \item[-] a transitive 
 set $\Delta_g$ homeomorphic to $\Lambda \times
 M$,
 \item[-] a homoclinic tangency (in $\Delta_g$) of codimension $c>0$.
 \end{enumerate}
The codimension $c$ of the tangency can be chosen to be any
integer $0<c\leq \dim M/2$.

\end{mainthm}

As far as we know, the theorem above gives the first direct
construction (i.e., not based on  a dimension reduction argument using
normally-hyperbolic manifolds) of robust tangencies for
symplectomorphisms in higher dimensions.

Notice that $f_0$ cannot be a complete integrable system since the
fiber map on $M$ is the identity function. Even so, we construct
nearby symplectomorphisms with a transitive set which projets onto
$M$. Large transitive sets with this property had been previously
constructed in~\cite{NP12} only close to integrable systems. The
integrability could be a restriction on the manifold $M$ since, as
far as we know, it is unknown if a given symplectic structure
admits integrable systems. Thus, Theorem~\ref{thmD} covers, a
priori, new examples where diffusion orbits can be obtained. To
prove this result we will use a different method to get orbits
drifting along $M$. The approach of~\cite{NP12} consists in
constructing transversal invariant tori for the fiber dynamics on
$M$. Also the geometrical mechanism of diffusion developed
in~\cite{DLS06} is based on a similar idea. Here we introduce a
different mechanism of propagation (drift) called
\emph{globalization}, which will be obtained by means of small
translations in Darboux local charts, compatible with the
symplectic structure.

The paper is organized as follows. First, in \S\ref{sec1} we
introduce blenders in a framework of normally hyperbolic invariant
laminations. These invariant laminations give rise to natural
skew-shifts over the space of sequences in a finite number of symbols,
called \emph{symbolic skew-products}. These kind of systems
have also been studied in~\cite{OGK17} and are related to the formation of
stochastic diffusive behaviour for the generalized Arnold example.
Next, in \S\ref{ss:construction-mixing} we provide a criterium to
construct robustly transitive symbolic skew-products. Finally, in
\S\ref{sec:sympletic} the main result is proven by combining the
criteria for robust transitivity in
\S\ref{ss:construction-mixing} and for robust tangencies
from~\cite{BR17}.

%% file: transitivity-preliminaries.tex
\section{Symbolic skew-products and blenders} \label{sec1}
Let  $\mathscr{A}$ be  a finite set (with at least two points),
that we call an alphabet of symbols, and fix $0<\nu<1$ and
$0<\alpha\leq 1$. Consider the product space
$\Sigma\equiv\Sigma(\mathscr{A},\nu)\eqdef \mathscr{A}^\mathbb{Z}$
of the bi-sequences $\xi=(\xi_i)_{i\in\mathbb{Z}}$ of symbols in
$\mathscr{A}$ endowed with the metric
\begin{equation*}
\label{e:metrica} d_{\Sigma}(\xi,\zeta)\eqdef \nu^{\ell}, \ \
\ell=\min\{i\geq 0: \xi_i\not=\zeta_i\ \text{or} \
\xi_{-i}\not=\zeta_{-i} \}.
\end{equation*}
In what follows $M$ will denote a differentiable manifold (not
necessarily compact and not necessarily boundaryless) of dimension
$c\geq 1$.

\subsection{Symbolic skew-products} \label{s:symbolic}
Given a compact set
$K$ in $M$, we consider the pseudometric in the set $C^0(M)$ of
continuous functions of $M$ given by
\begin{equation}
\label{eq:metrica00}
  d_{C^0}(\phi,\psi)_K \eqdef \max_{x,y\in K} d(\phi(x),\psi(x))
  \quad \text{for any $\phi, \psi \in C^0(M)$.}
\end{equation}
Since $M$ is $\sigma$-compact, there is a sequence of relatively
compact subsets $K_n$ whose union is $M$ and then we can endow
$C^0(M)$ with the weak topology (also called compact-open
topology) induced by the family of
pseudometrics~\eqref{eq:metrica00}. That is,
$$
   d_{C^0}(\phi,\psi)=\sum_{n=1}^\infty 2^{-n} \,
   \frac{d_{C^0}(\phi,\psi)_{K_n}}{1+d_{C^0}(\phi,\psi)_{K_n}},
   \quad
   \text{for any $\phi,\psi \in C^0(M)$}.
$$


\subsubsection{The set of  skew-products}

We consider skew-product homeomorphisms of the form
\begin{equation}
\label{e:sym-skew}
   \Phi: \Sigma \times M  \to \Sigma \times M, \qquad \Phi(\xi,x)=(\tau(\xi), \phi_\xi(x))
\end{equation}
where the base map $\tau:\Sigma\to\Sigma$ is the lateral shift map
and the fiber maps $\phi_\xi:M\to M$ are homeomorphisms of $M$. In
order to emphasize the role of the fiber maps we write
$\Phi=\tau\ltimes\phi_\xi$ and call it a \emph{symbolic
skew-product}. When no confusion arises we also write
$\mathcal{M}=\Sigma\times M$.

For every $n>0$ and $(\xi,x)\in \mathcal{M}$ set
\begin{equation*}
\label{n.seq}
\begin{aligned}
   \phi^n_\xi(x) \eqdef \phi_{\tau^{n-1}(\xi)}\circ\cdots\circ\phi_{\xi}(x) 
    \quad \text{and} \quad
    \phi^{-n}_\xi(x) \eqdef \phi^{-1}_{\tau^{-n}(\xi)}\circ\cdots
    \circ\phi^{-1}_{\tau^{-1}(\xi)}(x)
\end{aligned}
\end{equation*}
and hence
$$
\Phi^n(\xi,x)= (\tau^n(\xi),\phi^n_\xi(x)) \quad \text{for all
$n\in\mathbb{Z}$}.
$$

We introduce the set of symbolic skew-products with which we will
work:

\begin{defi}
Denote by
$\mathcal{S}(M)\equiv\mathcal{S}^\alpha_{\mathscr{A},\nu}(M)$ the
set of $\alpha$-H\"older continuous symbolic skew-products of
$\mathcal{M}=\Sigma\times M$. This is, the set of symbolic
skew-products $\Phi=\tau\ltimes \phi_\xi$ as in~\eqref{e:sym-skew}
such that
\begin{enumerate}[leftmargin=0.6cm, label=\textbullet]
\item
$\phi_\xi$ are  bi-Lipschitz homeomorphisms (uniform in $\xi$):
there are positive constants $\gamma\equiv\gamma(\Phi)>0$ and
$\hat\gamma\equiv \hat\gamma(\Phi)>0$ such that
\begin{equation}
\label{eq:lipschizt}
 \gamma\,d(x,y)< d(\phi_\xi(x),\phi_\xi(y)) <
\hat\gamma^{-1} \, d(x,y), \quad \mbox{ for all $x, y \in M$ and
$\xi \in \Sigma$,}
\end{equation}
\item $\phi_\xi$ depend $\alpha$-H\"older with respect to
$\xi$: there is a non-negative constant  $C_0 \equiv C_0(\Phi)\geq
0$ such that
\begin{equation} \label{eq:Holder} d_{C^0}(\phi^{\pm 1}_\xi,\phi^{\pm
1}_\zeta) \leq C_0 \, d_{\Sigma}(\xi,\zeta)^{\alpha} \quad
\mbox{for all $ \xi, \zeta \in \Sigma$ with $\xi_0=\zeta_0$.}
\end{equation}
\end{enumerate}
We will denote by $\mathcal{S}^{0}(M)$ the set $\mathcal{S}(M)$
with $C^1$-diffeomorphisms for fiber maps.
\end{defi}

We define in $\mathcal{S}(M)$ the metric
\begin{equation}
\label{eq:metrica0}
  d_{\mathcal{S}}(\Phi,\Psi) \eqdef
  d_0(\Phi,\Psi) + \mathrm{Lip}_0(\Phi,\Psi)+
  \mathrm{Hol}_0(\Phi,\Psi)
\end{equation}
where the symbolic skew-products $\Phi=\tau\ltimes\phi_\xi$ and
$\Psi=\tau\ltimes\psi_\xi$ belong to $\mathcal{S}(M)$ and
\begin{gather*}
\mathrm{Lip}_0(\Phi,\Psi)\eqdef \max_{\xi \in \Sigma}
  \big|\mathrm{Lip}(\phi_\xi^{\pm 1})-\mathrm{Lip}(\psi_\xi^{\pm
  1})\big|. \\
   d_0(\Phi, \Psi)\eqdef \max_{\xi \in \Sigma} \,
 d_{C^0}(\phi^{\pm 1}_\xi,\psi^{\pm 1}_\xi)
\quad  \text{and} \quad \mathrm{Hol}_0(\Phi,\Psi)\eqdef
\big|C_0(\Phi)-C_0(\Psi)\big|
\end{gather*}
with
$$
\mathrm{Lip}(\phi)=\sup_{ x \not = y}
\frac{d(\phi(x),\phi(y))}{d(x,y)} \quad \text{being $\phi$ a
bi-Lipschitz homeomorphism of $M$.}
$$

An important class of $\alpha$-H\"older continuous symbolic
skew-products is the following:

\begin{defi}
\label{def:PHS0} A symbolic skew-product $\Phi=\tau\ltimes
\phi_\xi \in \mathcal{S}(M)$ is \emph{partially hyperbolic} if
$$\nu^\alpha < \gamma < 1 < \hat{\gamma}^{-1}< \nu^\alpha$$
where $\gamma$ and $\hat\gamma$ are given in~\eqref{eq:lipschizt}.
 We
denote by
$\mathcal{PHS}(M)\equiv\mathcal{PHS}^\alpha_{\mathscr{A},\nu}(M)$
the set of partially hyperbolic symbolic skew-products. In
addition, $\mathcal{PHS}^{0}(M)=\mathcal{PHS}(M)\cap
\mathcal{S}^{0}(M)$.
\end{defi}

Finally we introduce a particular and important class of symbolic
skew-products:

\begin{defi}[one-step maps]
A symbolic skew-product $\Phi=\tau\ltimes\phi_\xi$ is called
\emph{one-step} if the fiber maps $\phi_\xi$ only depend on the
coordinate $\xi_0$ of the bi-sequences
$\xi=(\xi_i)_{i\in\mathbb{Z}} \in \Sigma$. In this case we have
$\phi_\xi=\phi_i$ if $\xi_0=i$ and write $\Phi=\tau\ltimes
(\phi_1,\dots,\phi_d)$.
\end{defi}

The \emph{iterated function system} (IFS for
short) associated with a one-step map
$\Phi=\tau\ltimes(\phi_1,\dots,\phi_d)$ is the semigroup action
generated by $\phi_1,\dots,\phi_d$. In what follows, $\langle
\phi_1,\dots,\phi_d\rangle^+$ will denote the semigroup generated
by the maps $\phi_1,\dots,\phi_d$.


\subsubsection{Stable and unstable sets for skew-products}

We define the \emph{local stable} and \emph{unstable} set of the
lateral shift map $\tau:\Sigma \to \Sigma$ at $\xi\in \Sigma$
respectively as
\begin{align*}
W^s_{loc}(\xi) \equiv W^s_{loc}(\xi;\tau) \eqdef \{\zeta \in
\Sigma : \zeta_i=\xi_i, \ i\geq 0 \}, \quad 
W^u_{loc}(\xi) \equiv W^u_{loc}(\xi;\tau) \eqdef \{\zeta \in
\Sigma: \zeta_i=\xi_i, \ i<0\}.
\end{align*}
The \emph{(global) stable set} of the skew-product map
$\Phi:\mathcal{M} \to \mathcal{M}$ at $P\in \mathcal{M}$ is
defined as
$$
W^{s}(P)\equiv W^{s}(P;\Phi)\eqdef \{Q\in \mathcal{M}:
\displaystyle\lim_{n\to\infty} d(\Phi^{n}(Q),\Phi^{n}(P)= 0 \}.
$$
We define the \emph{(global) stable set} of  a compact
$\Phi$-invariant set, i.e.~so that $\Phi(\Gamma)=\Gamma$,  by
$$
W^s(\Gamma) \equiv W^s(\Gamma;\Phi)\eqdef \{P\in \mathcal{M}:
    \lim_{n\to\infty} d(\Phi^{n}(P),\Gamma)= 0\}
$$
or equivalently as the set of the points of $\mathcal{M}$ so that
its $\omega$-limit is contained in $\Gamma$. The set $\Gamma$ is
called \emph{isolated} (or \emph{maximal invariant set}) if there
is a compact neighborhood $\mathcal{U}$ of $\Gamma$, called the
\emph{isolating neighborhood} for $\Gamma$, such that every
invariant subset of $\mathcal{U}$ lies in $\Gamma$. In such a
case, we introduce the \emph{local stable set} of $\Gamma$ as the
forward invariant set of $\Phi$ in the isolating neighborhood
$\mathcal{U}$, that is,
\begin{equation*}
W^s_{loc}(\Gamma) \equiv W^{s}_{loc}(\Gamma;\Phi) \eqdef
  \{P\in\mathcal{M}: \ \Phi^n(P)\in \mathcal{U}  \  \text{for $n\geq 0$}\} = \bigcap_{n\geq 0} \Phi^n(\mathcal{U}).
\end{equation*}
Similarly  $W^u_{loc}(\Gamma) \equiv W^u_{loc}(\Gamma;\Phi)$ and
$W^u(\Gamma) \equiv W^u(\Gamma;\Phi)$ are, respectively, the
\emph{local unstable set} and the \emph{global unstable set} of
$\Gamma$. We have that
$$
   W^s(\Gamma)=\bigcup_{n\geq 0} \Phi^{-n}(W^s_{loc}(\Gamma)) \quad \text{and} \quad W^u(\Gamma)= \bigcup_{n\geq0} \Phi^{n}(W^u_{loc}(\Gamma)).
$$

Finally, given an \emph{$\mathcal{S}$-perturbation} of $\Phi$,
that is a symbolic skew-product $\Psi$ close to $\Phi$ in the
metric given in~\eqref{eq:metrica0}, we denote by $\Gamma_\Psi$
the maximal invariant set in $U$ of $\Psi$. Although isolated sets
vary, a priori, just upper semicontinuously by an abuse of
terminology, we call $\Gamma_\Psi$ the \emph{continuation} of
$\Gamma$ for $\Psi$.

\subsubsection{Strong laminations for partially hyperbolic skew-products}
Under the global assumption of domination introduced in
Definition~\ref{def:PHS0}, the usual graph transform argument
yields a local strong stable $\mathcal{W}^{ss}$ and unstable
$\mathcal{W}^{uu}$ partitions:

\begin{prop}[\cite{AV10,ASV12}]
\label{p:inv-s} For every $\Phi\in\mathcal{PHS}(M)$ there exist
unique partitions
\begin{align*}
\mathcal{W}^{ss}=\{W^{ss}_{loc}(\xi,x): \, (\xi,x) \in \mathcal{M}
\} \quad \text{and} \quad \mathcal{W}^{uu}=\{W^{uu}_{loc}(\xi,x):
\, (\xi,x) \in \mathcal{M} \}
\end{align*}
of $\mathcal{M}=\Sigma \times M$ such that
\begin{enumerate}[itemsep=0.1cm]
\item \label{item-kk1} every leaf $W^{ss}_{loc}(\xi,x)$
 is the graph of an $\alpha$-H\"older function
 $\gamma^{s}_{\xi,x} : W^s_{loc}(\xi) \to M$,
\item $W^{ss}_{loc}(\xi,x)$
varies continuously with respect to $(\xi,x)$ and depends
continuously on $\Phi$,
\item \label{item-kk2} $\Phi(W^{ss}_{loc}(\xi,x)) \subset
W^{ss}_{loc}(\Phi(\xi,x))$ for all $(\xi,x) \in \mathcal{M}$,
\item $W^{ss}_{loc}(\xi,x) \subset
W^{s}(\xi,x)$ for all $(\xi,x) \in \mathcal{M}$.
\end{enumerate}
The partition $\mathcal{W}^{uu}$ verifies analogous properties.
\end{prop}

Each leaf of the partition $\mathcal{W}^{ss}$  is called the
\emph{local strong stable set}. We define the \emph{(global)
strong stable set} of $\Phi$ at $P$ as
\begin{equation*}
 W^{ss}(P) \equiv W^{ss} (P;\Phi ) \eqdef
\displaystyle\bigcup_{n\geq 0} \Phi^{-n} (
W^{ss}_{loc}(\Phi^{n}(P))) \subset W^s(P).
\end{equation*}

\subsection{Blenders}
\label{s:symbolic-skew-products}  In this subsection, we will
first introduce the notion of hyperbolic set for symbolic
skew-products homeomorphisms. After that we give the formal
definition of blenders and finally we provide a criterion to
obtain these local tools.

\subsubsection{Hyperbolic sets}
\label{ss:hyperbolic}
Fix $\varepsilon>0$ small enough. We introduce the \emph{local
stable set (of size $\varepsilon$)}  of $\Phi$ at $P=(\xi,x)$ as
$$
   W^{s}_{\varepsilon}(P)\equiv W^{s}_{\varepsilon}(P;  \Phi) \eqdef \{Q\in \mathcal{M}: \
   d(\Phi^n(Q),\Phi^n(P))\leq \varepsilon, \ \ n\geq 0
   \} \subset W^s_{loc}(\xi)\times M.
$$
The local unstable set (of size $\varepsilon$), denoted by
$W^u_\varepsilon(P)$, is defined analogously.

\begin{defi}
\label{def:hyp} A compact invariant set $\Gamma\subset
\mathcal{M}$ is \emph{hyperbolic} (for $\Phi$) if there exist
constants $\varepsilon>0$, $K>0$, $0<\theta<1$ such that
\begin{align*}
d(\Phi^n(P),\Phi^n(Q))\leq K \theta^n \ \ \text{for all
$P\in \Gamma$,  $Q\in W^s_{\varepsilon}(P)$ and $n\geq 0$;} \\
d(\Phi^{-n}(P),\Phi^{-n}(Q))\leq K \theta^n \ \ \text{for all
$P\in \Gamma$,  $Q\in W^u_{\varepsilon}(P)$ and $n\geq 0$;}
\end{align*}
and there exists $\delta>0$ such that
$$
\# \, W^s_{\varepsilon}(P) \cap W^u_{\varepsilon}(Q) =1 \ \
\text{for all $P, \, Q \in \Gamma$ with $d(P,Q) \leq \delta$}.
$$
\end{defi}

Every isolated hyperbolic set $\Gamma$ for $\Phi$ is topologically
stable~\cite{Akin93}; i.e., there is an isolating neighborhood $U$
of $\Gamma$ such that for any homeomorphism $\Psi$ which is $C^0$
near $\Phi$, the restriction of $\Psi$ to the maximal invariant
set in $U$, is semiconjugate to the restriction of $\Phi$ to
$\Gamma$.

We will now introduce the notion of index of an isolated
transitive hyperbolic set $\Gamma$ in our context. In the sequel
we will assume that the topological dimension (in the sense of the
Lebesgue covering dimension) of
$M^{cs}_{\varepsilon}(P)=W^{s}_\varepsilon(P)\cap (\{\xi\}\times
M)$ and $M^{cu}_{\varepsilon}(P)=W^{u}_\varepsilon(P)\cap
(\{\xi\}\times M)$ depend continuously with respect to
$P=(\xi,x)\in\Gamma$. From this assumption and being $\Gamma$
transitive, the dimensions of $M^{cs}_{\varepsilon}(P)$ and
$M^{cu}_\varepsilon(P)$ remain constant for any $P\in\Gamma$.
Thus, we may define the \emph{$cs$-index} and \emph{$cu$-index} of
$\Gamma$, denoted by $\mathrm{ind}^{cs}(\Gamma)$ and
$\mathrm{ind}^{cu}(\Gamma)$ as these dimensions respectively.
Notice that $\dim M=
\mathrm{ind}^{cs}(\Gamma)+\mathrm{ind}^{cu}(\Gamma)$ and from the
topological stability, the $cs$-index remains constant under small
$\mathcal{S}$-perturbations of $\Phi$.

\subsubsection{Symbolic blenders}
In order to introduce the notion of a blender we need first to
define families of $s$-discs and $u$-discs.
To do this, we will consider a \emph{basic open} set $\mathcal{B}$
of $\mathcal{M}$, i.e., a set of the form $\mathsf{V}\times B$
where $\mathsf{V}$ is an open set of $\Sigma$ and $B$ is an open
set of $M$.

\begin{defi}[$s$-discs]
\label{discs} A set $\mathcal{D}^s \subset \mathcal{M}$ is called
a \emph{$s$-disc} in $\mathcal{B}$ if there is  $\xi \in
\mathsf{V}$ such that $\mathcal{D}^s$ is a graph of an
$\alpha$-H\"older function from $W^s_{loc}(\xi)\cap \mathsf{V}$ to
$B$.
\end{defi}

We say that two $s$-discs, $\mathcal{D}^s_1, \mathcal{D}^s_2
\subset W^s_{loc}(\xi)\times M$ are close if they are the graphs
of close $\alpha$-H\"older functions. This proximity between discs
allows us to introduce the following:

\begin{defi}[open set of $s$-discs]
We say that a collection of discs $\mathscr{D}^s$ is an \emph{open
set of $s$-discs in $\mathcal{B}$} if given $\mathcal{D}^s_0\in
\mathscr{D}^s$, every $s$-disc $\mathcal{D}^s$ close enough to
$\mathcal{D}^s_0$ is a $s$-disc contained in $\mathcal{B}$ and
belongs~to~$\mathscr{D}^s$.
\end{defi}

Example of $s$-discs are the \emph{almost horizontal discs}
defined as follows: given $\delta>0$ and a point $(\xi, x)\in
\mathcal{B}$, we say that a set $\mathcal{D}^s\equiv
\mathcal{D}^s(\xi, x) \subset \mathcal{M}$ is a
\emph{$\delta$-horizontal disc} in $\mathcal{B}$ if
\begin{enumerate}
 \item[-] $\mathcal{D}^s$ is a graph of a $(\alpha,C)$-H\"older function
 $g : W^s_{loc}(\xi)\cap \mathsf{V} \to B$,
 \item[-] $d(g(\zeta),x)<\delta$ for all
 $\zeta\in W^s_{loc}(\xi)\cap \mathsf{V}$,
 \item[-] $C\nu^\alpha <\delta$.
\end{enumerate}
The set of all $\delta$-horizontal discs in $\mathcal{B}$ is an
open set of $s$-discs in $\mathcal{B}$. Similarly we define
\emph{$u$-discs} in $\mathcal{B}$, \emph{open set of $u$-discs} in
$\mathcal{B}$ and we have that the set of \emph{almost vertical
discs} is an example of an open set of $u$-discs.

Following~\cite{NP12,BKR14,BR17}, we introduce symbolic $cs$, $cu$
and $double$-blenders.

\begin{defi}[blenders] \label{d:symbolic-blender}
Let $\Phi \in\mathcal{S}(M)$ be a symbolic skew-product. A
transitive hyperbolic maximal invariant set $\Gamma$ in a
relatively compact open set $\mathcal{U}\subset \mathcal{M}=
\Sigma\times M$ of $\Phi$ is called
\begin{enumerate}
 \item \label{eq:cs-blender} $cs$-\emph{blender} if \,$\mathrm{ind}^{cs}(\Gamma)>0$ and there
 exist a basic open set $\mathcal{B} \subset \mathcal{U}$ and
 an open set $\mathscr{D}^s$ of $s$-discs in $\mathcal{B}$
 such that for every small enough
 $\mathcal{S}$-perturbation $\Psi$ of $\Phi$,
 \begin{equation*}
 W^{u}_{loc}(\Gamma_\Psi)\cap \mathcal{D}^s\neq\emptyset \quad \text{for all $\mathcal{D}^s \in \mathscr{D}^s$.}
 \end{equation*}
 \item  \label{eq:cu-blender} $cu$-\emph{blender} if \,$\mathrm{ind}^{cu}(\Gamma)>0$ and there
 exist a basic open set $\mathcal{B}\subset \mathcal{U}$ and
 an open set $\mathscr{D}^u$ of $u$-discs in $\mathcal{B}$
 such that for every small enough
 $\mathcal{S}$-perturbation $\Psi$ of $\Phi$,
 \begin{equation*}
 W^{s}_{loc}(\Gamma_\Psi)\cap \mathcal{D}^u\neq\emptyset \quad
 \text{for all $\mathcal{D}^u\in \mathscr{D}^u$}.
 \end{equation*}
 \item \emph{double-blender} if both~\eqref{eq:cs-blender} and~\eqref{eq:cu-blender} hold (not necessarily for the same $\mathcal{B}$).
\end{enumerate}
The open set $\mathcal{B}$ is called a \emph{superposition domain}
and the open sets of discs $\mathscr{D}^s$ and $\mathscr{D}^u$ are
called the \emph{superposition regions} of the blender. Finally,
the $cs$-blender (resp.~$cu$-blender) with $cs$-index
(resp.~$cu$-index) is equal to $\dim M$ is called a
$contracting$-blender (resp.~$expanding$-blender).
\end{defi}

Blenders are actually a power tool in partially hyperbolic
dynamics when the superposition region contains the local strong
stable/unstable set in the superposition domain. For this reason,
without loss of generality, we will assume the following blender
properties:

\begin{scho} Consider $\Phi \in
\mathcal{PHS}(M)$ and let $\Gamma$ be a $cs$-blender with
superposition domain $\mathcal{B}= \mathsf{V}\times B$ and
superposition region $\mathscr{D}^s$. Then, there is a
$\mathcal{S}$-neighborhood $\mathscr{U}$ of $\Phi$ such that for
any $\Psi\in \mathscr{U}$,
\begin{enumerate}[label=(B\arabic*),ref=B\arabic*]
\item \label{cond:B1}
\emph{the open set of discs $\mathscr{D}^s$ contains the family of
local strong stable sets of $\Psi$ in $\mathcal{B}$:
$$
\text{if \ $W^{ss}_{loc}(P;\Phi)\cap (\mathsf{V}\times M) \subset
\mathcal{B}$ \ \ then \ \ $W^{ss}_{loc}(P;\Psi) \cap
(\mathsf{V}\times M) \in \mathscr{D}^s$.}
$$}
\item \label{cond:B2}
\emph{if $W^{ss}_{loc}(P;\Phi)\cap (\mathsf{V}\times M) \subset
\mathcal{B}$ then
$$
         W^u_{loc}(\Gamma_\Psi;\Psi) \cap W^{ss}_{loc}(P';\Psi)
         \not=\emptyset \quad \text{for all $P'$ close enough to $P$.}
$$}
\end{enumerate}
Similar conditions are also assumed for $cu$-blenders of partially
hyperbolic skew-products.
\end{scho}

We must show that property~\eqref{cond:B2} follows from
the definition of a blender.
\begin{proof}
 First of all, notice that the
assumption~\eqref{cond:B1} and Definition~\ref{d:symbolic-blender}
imply that
\begin{equation}
\label{eq:3} \text{if $W^{ss}_{loc}(P;\Phi)\cap (\mathsf{V}\times
M) \subset \mathcal{B}$  \ then \
         $W^u_{loc}(\Gamma;\Phi) \cap W^{ss}_{loc}(P;\Phi)
         \not=\emptyset$ \quad $\mathcal{S}$-robustly.
         }
\end{equation}
A priori, the neighborhood of the $\mathcal{S}$-perturbation of
$\Phi$ where~\eqref{eq:3} holds depends on the $s$-disc
$W^{ss}_{loc}(P;\Phi)$. However, this can be taken independent of
the disc assuming that the disc belongs to a superposition
subdomain $\mathcal{B}_0$ of $\mathcal{B}=\mathsf{V}\times B$.
That is if $ W^{ss}_{loc}(P;\Phi)\cap (V\times M) \subset
\mathcal{B}_0$ where $\mathcal{B}_0=V\times B_0$, $B_0$ is an open
set whose closure is contained in $B$. For this reason, without
loss of generality, we can assume that~\eqref{cond:B2} holds.
\end{proof}

\subsubsection{Blenders from one-step maps}
In order to provide a criterion to construct blenders we need the
following definition.

\begin{defi}[blending region]
\label{def:blending-region1} Let $\Phi=\tau\ltimes
(\phi_1,\dots,\phi_d)\in \mathcal{PHS}(M)$. Consider  bounded open
sets $B$ and $D$ of $M$ with $\overline{B}\subset D$, a subset
$S\subset \mathscr{A}$, and a hyperbolic transitive set
$$
\Gamma \eqdef \bigcap_{n\in\mathbb{Z}} \Phi^n(S^\mathbb{Z}\times
\overline{D})=\bigcap_{n\in\mathbb{Z}} \Phi^n(\mathsf{V}\times
\overline{D}),
$$
where $\mathsf{V}$ denotes any isolating neighborhood of
$\Sigma^+_S\eqdef\{\xi\in \Sigma: \xi_0 \in S\}$ and
 $\Sigma^-_S\eqdef\{\xi\in \Sigma: \xi_{-1} \in S\}$.
We say that $B$ is a \emph{$cs/cu/double$-blending region} with
respect to $\{ \phi_i: i\in S\}$ on $D$  if there exists
respectively a
\begin{enumerate}
\item \label{cs:cover} \emph{$cs$-cover}: $\{\phi_i(B):i \in S \}$ is an open cover of $\overline{B}$ and $\mathrm{ind}^{cs}(\Gamma)>0$;
\item \label{cu:cover} \emph{$cu$-cover}: $\{\phi_i^{-1}(B): i \in S\}$ is an open cover of $\overline{B}$ and $\mathrm{ind}^{cu}(\Gamma)>0$;
\item \emph{$double$-cover}: both~\eqref{cs:cover} and~\eqref{cu:cover} are true.
\end{enumerate}
We call \emph{$cs$-index} (resp.~\emph{$cu$-index}) of the
blending region $B$ the $cs$-index (resp.~$cu$-index) of $\Gamma$.
As in the case of the blender, if its $cs$-index
(resp.~$cu$-index) is equal to dimension of $M$ the blending
region is called \emph{contracting} (resp.~\emph{expanding}).
\end{defi}

With the above terminology, the following result showed
in~\cite[Corolory~5.3]{BR17} gives a criterion to construct a
blender.
\begin{thm}\label{cor:blender} Let $\Phi=\tau\ltimes(\phi_1,\dots,\phi_d) \in
\mathcal{PHS}(M)$. Assume that there are a set $S\subset
\mathscr{A}$ and a
\begin{itemize}
 \item[-]  $cs/cu/double$-blending region $B$ with respect to $\{ \phi_i: i\in S\}$ on $D$.
\end{itemize}
Then the maximal invariant set $\Gamma$ in $S^\mathbb{Z}\times
\overline{D}$ is a $cs/cu/double$-blender  of $\Phi$ whose
super\-po\-si\-tion region contains the family of almost
horizontal discs in $\Sigma^+_S\times B$ or/and almost vertical
discs in $\Sigma^-_S\times B$. Moreover, it also contains the
family of local strong stable/unstable sets, i.e.,~\eqref{cond:B1}
holds.
\end{thm}

Blending regions which cover an open ball
$B$ around a hyperbolic fixed point of a map $\phi$ can be easily constructed from a
sufficient number of sets of the form $\phi_i(B)$, where $\phi_i$
is a translation of $\phi$. This idea was developed
in~\cite{NP12} and ~\cite[Proposition~5.6]{BR17}, obtaining a blending region in
local coordinates:

\begin{prop}
\label{prop:creation-blender} Consider a $C^r$-diffeomorphism
$\phi$ of $\mathbb{R}^c$ with a hyperbolic
attrac\-ting/repe\-lling/saddle fixed point $x$. Then, there exist
an integer $k\equiv k(\phi,c)\geq 2$, arcs of
$C^r$-diffeomorphisms of $\mathbb{R}^c$,
$\phi_1\equiv\phi_1(\varepsilon),\dots,\phi_k\equiv\phi_k(\varepsilon)$
and bounded open sets $D\equiv
D(\varepsilon)$,  $\varepsilon\geq 0$, 
such that
\begin{enumerate}[leftmargin=0.6cm]
\item[-] $\phi_i(0)=\phi$ for $i=1,\dots,k$;
\item[-] $\phi_i=T_i \circ \phi$ where $T_i\equiv T_i(\varepsilon)$ is a translation (moreover, one can take $\phi_1=\phi$);
\item[-] $B\equiv B_{\delta}(x) \subset D \subset B_{2\varepsilon}(x)$ for some $\delta\equiv \delta(\varepsilon)>0$; 
\item[-] $B$ is $cs/cu/double$-blending region with respect to $\{\phi_1,\dots,\phi_k\}$ on $D$ for all $\varepsilon>0$.
\end{enumerate}
Moreover, the $cs$-index of the blending region is equal to the
$s$-index of the hyperbolic fixed point $x$.
\end{prop}

%% file: transitivity-criteria.tex
\section{Robust transitivity}
\label{ss:construction-mixing} We explain how blenders can be used
to yield a \emph{$\mathcal{S}^0$-robust topologically mixing}
symbolic skew-product $\Phi$ of $\mathcal{M}=\Sigma\times M$. That
is, for any $\mathcal{S}^0$-perturbation $\Psi$ of $\Phi$ and for
every pair of open sets $\mathcal{U}$, $\mathcal{V}$ of
$\mathcal{M}$, there is $n_0> 0$ such that $\Psi^n(\mathcal{U})
\cap \mathcal{V} \not=\emptyset$ for all $n \geq n_0$. In
particular topologically mixing implies transitivity.

\subsection{Criterion to yield robust transitivity} One of the classical ways to create robustly transitive
diffeomorphisms is to construct a map that robustly has a
hyperbolic periodic point with dense stable and unstable
manifolds. Then, using the inclination lemma (or $\lambda$-lemma)
one concludes that the diffeomorphism is topologically mixing. In
the symbolic setting an analogous result was proved
in~\cite[Theorem~5.7]{BKR14} for fiber attracting/repelling
hyperbolic fixed points. Here we extend this criterion to any
hyperbolic periodic point. Later on we will use this result to
construct robustly topologically mixing skew-products.

\begin{thm}[criterion for robust transitivity]
\label{thm:top-mixing} Assume that $\Phi \in \mathcal{PHS}^0(M)$
has a $cs$-blender $\Gamma$ with superposition domain
$\mathcal{B}=\mathsf{V}\times B$ such that
\begin{enumerate}[label=(RT),ref=RT]
 \item \label{RT1}
there is an open set $\mathcal{B}_0\subset \mathcal{B}$ such that
 $$
 \text{$\mathcal{B}_0=\mathsf{V}\times B_0$ with
 $\overline{B_0} \subset B$}, \quad
     W^{ss}_{loc}(Q)\cap (\mathsf{V}\times M) \subset \mathcal{B}
     \quad \text{for all $Q\in \mathcal{B}_0$}
 $$
 and
\begin{equation}
\label{RT}
    \mathcal{M} = \overline{\bigcup_{n=1}^\infty \Phi^n(\mathcal{B}_0)}  \qquad \text{$\mathcal{S}^0$-robustly.\footnotemark}
\footnotetext{
    For every compact set $K \subset M$ there is
    a $\mathcal{S}^0$-neighborhood $\mathscr{U}$ of $\Phi$ such that
    $\Sigma\times K$ is contained in the closure of
    $\cup_n \Psi^n(\mathcal{B}_0)$ for all
    $\Psi \in \mathscr{U}$.}
\end{equation}
\end{enumerate}

\noindent Then  $\mathcal{S}^0$-robustly it holds that
$\mathcal{M} = \overline{W^u(P)}$ for any periodic point $P\in
\Gamma$.
If $\Gamma$ is a $cu$-blender satisfying~\eqref{RT1} for
$\Phi^{-1}$, then the stable set of any periodic point of $\Gamma$
is $\mathcal{S}^0$-robustly dense in $\mathcal{M}$.

Moreover,  if  $\Gamma$ is either a $double$-blender, a
$contracting$-blender or an  $expanding$-blender
satisfying~\eqref{RT1} for both $\Phi$ and $\Phi^{-1}$, then the
global stable and unstable sets of any periodic point of $\Gamma$
are both $\mathcal{S}^0$-robustly dense. In this case, some power
$\Phi^k$ is $\mathcal{S}^0$-robustly topologically mixing which in
particular implies that $\Phi$ is $\mathcal{S}^0$-robustly
transitive.
\end{thm}

Notice that the globalization property~\eqref{RT} is a necessary
condition for robust transitivity. In order to prove
Theorem~\ref{thm:top-mixing} we need the following lemmas.

\begin{lem}
\label{lem:per-points} Consider a symbolic skew-product $\Phi$ and
let $\Gamma$ be an isolated hyperbolic transitive set of
$\Phi^\ell$ for some $\ell\in\mathbb{N}$.  Then for every periodic
point $P\in \Gamma$,
$$
    W^u(\Gamma;\Phi^\ell)\subset \overline{W^u(P;\Phi)} \quad
    \text{and}  \quad W^s(\Gamma;\Phi^\ell)\subset
    \overline{W^s(P;\Phi)}.
$$
\end{lem}
\begin{proof}
The transitivity of $\Gamma$ implies that $\Gamma \subset
\overline{W^u(P;\Phi^\ell)}$. Since $\Gamma$ is an isolated
hyperbolic set of $\Phi^\ell$, from the in-phase result
(c.f.~\cite[Prop.~10]{Ombach96}),
 it holds that
$$
   W^u(\Gamma;\Phi^\ell)=\bigcup_{Q\in\Gamma} W^u(Q;\Phi^\ell).
$$
Combining the above facts we get the density of the unstable set, and an
analogous argument concludes the density of the stable set.
\end{proof}

\begin{lem}
\label{lem:lambda-lema} Consider $\Phi\in\mathcal{S}^0(M)$ and let
$P$ be a hyperbolic periodic point of period $k>0$ such that
$W^u(P)\cap \mathcal{U} \not = \emptyset$ and $W^s(P)\cap
\mathcal{V}\not=\emptyset$ where $\mathcal{U}$ and $\mathcal{V}$
are open sets of $\mathcal{M}$. Then, there exists
$n_0\in\mathbb{N}$ such that
$$\Phi^{kn}(\mathcal{V})\cap \mathcal{U}
\not=\emptyset  \quad \text{for all $n\geq n_0$.}
$$
\end{lem}

\begin{proof}
Without loss of generality, we can assume that $P=(\xi,x)$ is a
fixed point of $\Phi$. Taking an integer $m>0$ large enough, we
get that $\Phi^{-m}(\mathcal{U})\supset
(W^{s}_{loc}(\xi)\cap \mathsf{C}) \times U$ and
$\Phi^{m}(\mathcal{V}) \supset (W^{u}_{loc}(\xi)\cap
\mathsf{C}) \times V$ where $\mathsf{C}$ is a cylinder around
$\xi$ and $U$, $V$ are open sets of $M$ so that
$$
  W^u_{\varepsilon}(x;\phi_{\xi}) \cap U \not=\emptyset
  \quad \text{and} \quad W^s_{\varepsilon}(x;\phi_{\xi}) \cap V \not=\emptyset
$$
for some $\varepsilon>0$ sufficiently small. Thus, it suffices to
prove the result for the hyperbolic fixed point $x\in M$  of the
$C^1$-diffeomorphism $\phi_\xi:M\to M$. Hence, as a consequence of
the inclination lemma we get that $\phi_{\xi}^{-n}(U) \cap
V\not=\emptyset$ for any $n>0$ large enough.
\end{proof}

\begin{proof}[Proof of Theorem~\ref{thm:top-mixing}]
Consider any open set $\mathcal{U} \subset \mathcal{M}$.
By~\eqref{RT}, there is $n$ arbitrarily large and $R\in
\mathcal{U}$ such that $Q=\Phi^{-n}(R) \in \mathcal{B}_0$. Hence,
$W^{ss}_{loc}(Q) \subset \Phi^{-n}(W^{ss}_{loc}(R)\cap
\mathcal{U})$. From~\eqref{RT1}, it holds $W^{ss}(Q) \cap
(\mathsf{V}\times M) \subset \mathcal{B}$. Since $\Gamma$ is a
$cs$-blender with superposition domain $\mathcal{B}$,
by~\eqref{cond:B2},
\begin{equation}
\label{eq:inter} W^{u}_{loc}(\Gamma)\cap
W^{ss}_{loc}(Q)\neq\emptyset \qquad
\text{$\mathcal{S}^0$-robustly.}
\end{equation}
Hence, Lemma~\ref{lem:per-points} implies $ W^{u}(P)\cap
\Phi^{-n}(\mathcal{U}) \neq \emptyset$ for all periodic points
$P\in \Gamma$. Consequently, $\mathcal{S}^0$-robustly it holds
$W^{u}(P)\cap \mathcal{U} \not \neq \emptyset$ . From the
invariance of the local strong stable partition, we get that
$W^{u}(P)$ meets $\mathcal{U}$, $\mathcal{S}^0$-robustly.

A similar argument holds assuming that $\Gamma$ is  a $cu$-blender
and~\eqref{RT1} holds for $\Phi^{-1}$. Thus according to
Lemma~\ref{lem:lambda-lema} we conclude the theorem when $\Gamma$
is a double-blender. If $\Gamma$ is a $contracting$-blender
(resp.~$expanding$-blender), then the local stable
(resp.~unstable) set of any periodic point in $\Gamma$ has
dimension $\dim M$. Hence, any $u$-disc (resp.~$s$-disc) in the
superposition region transversally meets the local stable
(resp.~unstable) set of these periodic points and
thus~\eqref{eq:inter} immediately holds for $\Phi^{-1}$ (resp.~for
$\Phi$). Therefore, in this case, the same argument works assuming
the globalization property~\eqref{RT1} for $\Phi$ and $\Phi^{-1}$.
This shows the density of both $W^s(P)$ and $W^u(P)$ and again,
according to Lemma~\ref{lem:lambda-lema}, we conclude the theorem.
\end{proof}

\begin{rem}
\label{rem:periodo} We actually get that $\Phi$ is
$\mathcal{S}^0$-robustly topologically mixing if $\Gamma$ has
fixed points.
\end{rem}




\subsection{Robust transitivity from one-step maps}
Theorem~\ref{thm:top-mixing} provides conditions to yield robustly
transitive symbolic skew-products. In what follows, we will first
translate these conditions to the particular case of one-step
maps, $\Phi=\tau\ltimes(\phi_1,\dots,\phi_d)$.  Afterwards, we will
construct arcs of IFSs unfolding from the identity and satisfying
these conditions.

\subsubsection{Criterion for robustly transitive one-step maps}
Let $\Phi=\tau\ltimes(\phi_1,\dots,\phi_d) \in \mathcal{PHS}(M)$
be a one-step map and assume there exists a $cs$-blending
region $B$ with respect to $\{\phi_1,\dots,\phi_d\}$ on $D$. Hence,
the maximal invariant set $\Gamma$ in the closure of $\Sigma\times
D$ for $\Phi$ is a $cs$-blender with superposition domain
$\mathcal{B}=\Sigma \times B$.

It is not difficult to see that~\eqref{RT} is equivalent to
\begin{equation}
\label{equiv-forward-globalization} M = \mathscr{P}\bigg(\
\overline{\bigcup_{n=1}^\infty \Phi^n(\mathcal{B}_0)}\ \bigg)
\quad \text{$\mathcal{S}^0$-robustly,}
\end{equation}
where $\mathcal{B}_0$ is an open set in $\mathcal{B}$
such that 
$$
\text{$\mathcal{B}_0=\Sigma\times B_0$ \ with \
$\overline{B_0}\subset B$,} \ \
     W^{ss}_{loc}(\xi,x) 
     \subset \mathcal{B} \ \ \text{for all $(\xi,x)\in \mathcal{B}_0$,
     \quad  $\mathcal{S}^0$-robustly}.
$$
The robustness in~\eqref{equiv-forward-globalization} means that
for every compact set $K \subset M$ there is a
$\mathcal{S}^0$-neighborhood $\mathscr{U}$ of $\Phi$ such that $K$
is contained in the projection on $M$ of the closure of the union
of forward $\Psi$-iterated of $\mathcal{B}_0$ for all $\Psi \in
\mathscr{U}$. Observe that~\eqref{equiv-forward-globalization}
holds if there exists $n\in\mathbb{N}$ such that
$$
      M=\mathscr{P}\big(\,\bigcup_{i=1}^{n} \Phi^{i}(\mathcal{B}_0) \,\big).
$$
This finite cover requires the compactness of $M$. In the
non-compact case, \eqref{equiv-forward-globalization} holds if
there exists an increasing sequence of compact sets $K_i\subset M$
such that their union is $M$ and for all $i$
$$
     K_i \subset \mathscr{P}\big( \, \bigcup_{n=1}^\infty \Phi^n(\mathcal{B}_0) \, \big)= \{h(x): h \in \langle\phi_1,\dots,\phi_d\rangle^+ \ \text{and} \ x\in B_0 \}
      \eqdef \langle \phi_1,\dots,\phi_d\rangle^+(B_0)
$$
where recall that $\langle \phi_1,\dots,\phi_d \rangle^+$ denotes
the semigroup generated by the maps $\phi_1,\dots,\phi_d$.

\begin{defi}[globalization]
Let $B$ be an open set of $M$. We say that $B$ is \emph{forward
and backward globalized} (or simply \emph{globalized})  by
$\langle \phi_1,\dots,\phi_d\rangle^+$ if for any compact set $K
\subset \mathrm{int}(M)$,
$$
  K \subset \langle \phi_1,\dots,\phi_d\rangle^+(B)
  \quad \text{and} \quad
  K \subset \langle \phi_1^{-1},\dots,\phi_d^{-1}\rangle^+(B).
$$
\end{defi}

With this terminology, as a consequence of
Theorem~\ref{thm:top-mixing} and Remark~\ref{rem:periodo} we have

\begin{cor}[robustly transitive one-step maps]
\label{cor:transitividad}
 Let $\Phi=\tau\ltimes
(\phi_1,\dots,\phi_d)\in \mathcal{PHS}^0(M)$. Suppose that there
are bounded open sets $D$, $B_0$ and $B$ of $M$  with
$\overline{B_0}\subset B\subset D$ such that
\begin{itemize}
\item[-] $B_0$ is globalized by $\langle
\phi_1,\dots,\phi_d\rangle^+$,
\item[-] $B$ is a $contracting/expanding/double$-blending
region with respect to $\{\phi_1,\dots,\phi_d\}$ on $D$.
\end{itemize}
Then $\Phi$ is $\mathcal{S}^0$-robustly transitive. Moreover,
$\Phi$ is $\mathcal{S}^0$-robustly topologically mixing if some
fiber map $\phi_i$ has a hyperbolic fixed point in $B$.
\end{cor}

\subsubsection{Construction of globalized blending regions}

Consider a $contrac\-ting/expan\-ding/dou\-ble$-blending region
$B$ with respect to $\{ \phi_1,\dots,\phi_d\}$ where $\phi_i$ are
diffeomorphisms of $M$. According to~\cite{BKR14,HN13} one way to
yield globalization is to add, if necessary, a pair of
Morse-Smale diffeomorphisms without periodic points in common and
having respectively an attracting/repelling periodic point in $B$
with a dense stable/unstable manifold. The existence of such a map
can be constructed by perturbations of time-one maps of
gradient-like vector fields.
This explains the following proposition.

\begin{prop}[globalized blending regions homotopic to the identity]
\label{thm:globalizacion-id} Let $M$ be a connected manifold of
dimension $c\geq 1$.  Then there exists an integer $s\equiv
s(c)\geq 3$ such that for every $x\in M$ there are arcs of
$C^r$-diffeomorphisms $T_1(\varepsilon), \dots, T_s(\varepsilon)$,
$\varepsilon\geq 0$, of $M$  such that
\begin{enumerate}[leftmargin=0.6cm]
\item[-] $T_i(0)=\mathrm{id}$  \ for $i=1,\dots,s$, and
\item[-] $B_{\varepsilon}(x)$ is globalized by $\langle T_1(\varepsilon), \dots, T_s(\varepsilon)\rangle^+$
 for all $\varepsilon>0$ small enough.
\end{enumerate}
\end{prop}


Next, we will give another proof of the above statement,
which will be useful later on for the symplectic setting. The following
construction uses \emph{local} tools in contrast to the
\emph{global} nature of Morse-Smale diffeomorphisms. The local
perturbations will be translations of the identity map,
compatible with symplectomorphisms. To show the result, we need
the following lemma and notation. Given $\delta>0$ and a subset
$A$ of $\mathbb{R}^c$, we write
$$
    B_\delta(A)\eqdef \{ y\in \mathbb{R}^c: d(y,A)<\delta \}.
$$
\begin{lem}
\label{lem:Rn} Let $U_0$ be a bounded connected open set of
$\mathbb{R}^c$. Then, there are an integer $m\equiv m(c)\geq 2$
(that only depends on the dimension $c$) and arcs of
$C^r$-diffeomorphisms $T_1(\varepsilon),\dots,T_m(\varepsilon)$ of
$U\equiv B_{2\varepsilon}(U_0)$ such that
for every $x\in U_0$, 
\begin{enumerate}[leftmargin=0.6cm, itemsep=0.1cm]
\item[-] $T_i(0)=\mathrm{id}$  \ for $i=1,\dots,m$,
\item[-] $T_i(\varepsilon)=\mathrm{id}$ in $\partial U$  for all
$i=1,\dots,m$ and $\varepsilon>0$,
\item[-] $T_i(\varepsilon)$ is a small translation in $U_0$
for all $i=1,\dots,m$ and $\varepsilon>0$,
\item[-] $\overline{U_0} \subset
\langle T_1(\varepsilon), \dots,
T_m(\varepsilon)\rangle^+(B_{\varepsilon}(x))$ for all
$\varepsilon>0$ small enough.
\end{enumerate}
\end{lem}
\vspace{-0.3cm}
\begin{proof}
Fix $\varepsilon>0$. Let $B\equiv B(\varepsilon)$ be an open ball
of radius $\varepsilon$ centered at a point $x$ in $U_0$. We can
cover the closure of this ball by small translations of $B$.
Notice that the number $m$ of translation required only depends on
the dimension $c$. In fact, there exist unitary vectors
$u_1,\dots,u_m$ so that \emph{any} ball in $\mathbb{R}^c$ can be
covered by translations along the directions $u_1,\dots,u_m$. Let
$T_{v_1},\dots,T_{v_m}$ be these translations of small vectors
$v_i=\delta u_i$, $0<\delta<\varepsilon$, such that $\overline B
\subset T_{v_1}(B) \cup \dots \cup T_{v_m}(B)$ with
$T_{v_i}(U_1)\subset B_{\varepsilon/2}(U_1)$ where $U_1\equiv
B_\varepsilon(U_0)$. By the above observation, the same
translations $T_{v_1},\dots, T_{v_m}$ cover any open ball of
radius $\rho \geq \varepsilon$.

By means of the extension lemma~\cite[Lemma~2.27]{Lee13}, there
exists a $C^r$-diffeomorphism $T_i\equiv T_i(\varepsilon)$ of $U$
such that $T_i|_{U_1}=T_{v_i}$ and $T_i|_{\partial U}=
\mathrm{id}$.
We will now show that $\overline{U_0} \subset \langle T_1, \dots,
T_m\rangle^+(B)$. Indeed, there is~$\rho>\varepsilon$~so~that
$$
     \overline{B_\rho}  \subset T_1(B)\cup\dots\cup T_m(B)
     \quad \text{where $B_\rho\equiv B_\rho(x)$}.
$$
If $B_\rho \subset U_1$, as the translation directions $u_i$ only
depend on the dimension $c$  then
$$
\overline{B_{2\rho}} \subset T_1(B_\rho)\cup \dots \cup
T_m(B_\rho).
$$
Repeating the above procedure, since the radius of the covered
ball $B_\rho$ is strictly increasing, by changing the center we
can reach the boundary of $U_1$ and thus cover $\partial U_0$ .

Since $v_i(\varepsilon)$ tends continuously to zero (i.e.
$\delta<\varepsilon$ goes to zero continuously), we get that
$T_i(\varepsilon)$ tends continuously to the identity and conclude
the proposition.
\end{proof}

\vspace{-0.3cm}
\begin{proof}[Proof of Proposition~\ref{thm:globalizacion-id}]
By means of the well-known procedure of Milnor
(see~\cite[Proposition~2.1]{Jam78}
and~\cite[Proposition~1.4.14]{MD92}), there is an atlas
$$\mathscr{A}=\{U_{ij}: i \in \mathbb{N}, \ j=0,\dots,\dim M\}$$ of
$M$ such that $\{U_{ij}\}_{i\in\mathbb{N}}$ are pairwise disjoint
for all $j=0,\dots, \dim M$. Let
$$
\mathscr{A}_0=\{U_{ij}^0: i \in \mathbb{N}, \ j=0,\dots,\dim M \}
\quad \text{with \ \ $U_{ij}=B_{\varepsilon_{ij}}(U_{ij}^0)$}$$ be
a refinement of $\mathscr{A}$. Relabeling the atlas if necessary,
we assume that $x\in U_{10}^0$. Lemma~\ref{lem:Rn} provides (in
local coordinates) $C^r$-diffeomorphisms
$T_{1}(\varepsilon_{10}),\dots,T_{m}(\varepsilon_{10})$ of
$U_{10}$, and an open ball $B(\varepsilon_{10})$ centered at $x$
such that $T_\ell(\varepsilon_{10})|_{\partial
U_{10}}=\mathrm{id}$ for all $\ell=1,\dots,m$ and
$$
\overline{U^0_{10}}\subset \langle
T_{1}(\varepsilon_{10}),\dots,T_{m}(\varepsilon_{10})\rangle^+(B(\varepsilon_{10})).
$$
For each $U_{ij}^0$ such that $U_{ij}^0\cap
U_{10}^0\not=\emptyset$, we take a point in this intersection and
apply again Lemma~\ref{lem:Rn} obtaining now new
$C^r$-diffeomorphisms
$T_{1}(\varepsilon_{ij}),\dots,T_{m}(\varepsilon_{ij})$ of
$U_{ij}$ and an open ball $B(\varepsilon_{ij}) \subset
U_{ij}^0\cap U_{10}^0$ such that
\begin{equation}
\label{eq:T}
T_\ell(\varepsilon_{ij})|_{\partial U_{ij}}=\mathrm{id} \ \
\text{for all $\ell=1,\dots,m$} \quad \text{and} \quad
\overline{U^0_{ij}}\subset \langle
T_{1}(\varepsilon_{ij}),\dots,T_{m}(\varepsilon_{ij})\rangle^+(B(\varepsilon_{ij})).
\end{equation}
Let $U_{i_1j_1}^0$ be one of the open sets choosen from the
preceeding step. Consider the open sets $U_{ij}^0$ such that
$U_{ij}^0\cap U_{i_1j_1}^0\not=\emptyset$ and which were not
considered in the previous step. Repeat the whole process by
choosing a point in the intersection and applying
Lemma~\ref{lem:Rn}. By induction, for all $i\in \mathbb{N}$ and
$j=0,\dots,\dim M$, there exist $C^r$-diffeomorphisms
$T_{1}(\varepsilon_{ij}),\dots,T_{m}(\varepsilon_{ij})$ of
$U_{ij}$ and an open ball $B(\varepsilon_{ij}) \subset U_{ij}^0$
satisfying~\eqref{eq:T}.

For each $\ell=1,\dots,m$ and $j=0,\dots,\dim M$, we define the
$C^r$-diffeomorphisms  $T_{\ell j}$ of $M$ by
$$
T_{\ell j}|_{U_{ij}^0}=T_{\ell}(\varepsilon_{ij}) \ \text{for all
$i\in \mathbb{N}$ and} \ \   T_{\ell j}=\mathrm{id} \ \text{in} \
M\setminus \bigcup_{i\in\mathbb{N}} U_{i j}.
$$
Observe that $T_{\ell j}$ is well defined since for each $j$,
$U_{ij}$, $i\in\mathbb{N}$ are pairwise disjoint open sets and
$T_{\ell}(\varepsilon_{ij})$ restricted to $\partial U_{ij}$ is
equal to the identity for all $\ell=1,\dots,m$. Moreover, $\langle
T_{\ell j} : \ \ell=1,\dots,m, \ j=0,\dots,\dim M \rangle^+$ has
forward globalization of the neighborhood $B(\varepsilon_{10})$ of
any $x$. Similarly by the same procedure, we can assume that this
semigroup
also has backward globalization of $B(\varepsilon_{10})$. Finally,
if $\varepsilon_{ij}\to 0$ then $T_{\ell j}$ tends continuously to
the identity  and thus relabeling the maps, we have arcs of
$C^r$-diffeomorphisms
$T_1=T_1(\varepsilon),\dots,T_s(\varepsilon)$, $\varepsilon\geq
0$, $s=(\dim M +1) m$, satisfying the required properties.
\end{proof}

\subsubsection{Arcs of robustly transitive one-step maps} The previous
proposition allows us to construct an arc of robustly transitive symbolic
skew-products.

\begin{thm}
\label{thm:mainA-symbolic-setting}
  Let $M$ be a connected manifold of dimension $c\geq 1$.  Then,
  there are an integer  $d\equiv d(c)\geq 3$ and an arc of
  one-step maps
  $\Phi_\epsilon=\tau\ltimes (\phi_1, \dots, \phi_d)$ in
$\mathcal{PHS}^r(M)$ isotopic to $\Phi_0=\tau\times \mathrm{id}$
so that $\Phi_\varepsilon$ is $\mathcal{S}^0$-robustly
topologically  mixing for any $\varepsilon>0$.
\end{thm}

\begin{proof}
Consider $x\in M$. By means of an arbitrarily small perturbation
of the identity map~\cite{H02,HN13} we can create a map $\phi$ for
which $x$ is a hyperbolic fixed point. Hence, applying
Proposition~\ref{prop:creation-blender} we can get arcs of
$C^r$-diffeomorphisms $\phi_1\equiv
\phi_1(\varepsilon),\dots,\phi_{k}\equiv \phi_{k}(\varepsilon)$
homotopic to the identity as $\varepsilon \to 0^+$, where $k\geq
2$ and only depends on $c$ and a $cs$-blending regions $B$ in
$B_{2\varepsilon}(x)$. Without loss of generality, assume that $B$
is a $contracting/double$-blending region. According to
Proposition~\ref{thm:globalizacion-id}, there exist $s\equiv
s(c)\geq 3$ arcs of $C^r$-diffeomorphisms $T_1\equiv
T_1(\varepsilon),\dots, T_s\equiv T_s(\varepsilon)$ homotopic to
the identity as $\varepsilon \to 0^+$ so that $\langle T_1,\dots,
T_s\rangle^+$ has globalization of a small open ball
$B_0\subset\overline{B}$. Adding these diffeomorphisms to the
previous maps if necessary, $B$ is a globalized
$contracting/double$-blending region. Take $d=s+k$, which only
depends on $c$, and set
$$\Phi_\varepsilon=\tau \ltimes
(\phi_1,\dots,\phi_{k},T_1,\dots,T_s) \in \mathcal{PHS}^r(M).$$
Hence $\Phi_0=\tau\times\mathrm{id}$ and according to
Theorem~\ref{cor:transitividad}, for any $\varepsilon>0$,
$\Phi_\varepsilon$ is $\mathcal{S}^0$-robustly topologically
mixing for any $\varepsilon>0$. This completes the proof of the
theorem.
\end{proof}

%
%
%

%% file: transitivity-symplectic.tex
\section{Symplectic skew-products}
\label{sec:sympletic} By a symplectic manifold $M$ we mean a
manifold equipped with a closed non-degenerate differential
two-form which is called the symplectic form. The nondegeneracy of
the form implies that the space must be even-dimensional. We will
create robust homoclinic tangencies inside semi-local
transitive partially hyperbolic sets for
\emph{symplectomorphisms}, that is diffeomorphisms preserving the
symplectic form. We will first show that the robust transitivity
of Theorem~\ref{thm:mainA-symbolic-setting} and robust tangencies
constructed in~\cite{BR17}
hold for \emph{symplectic
symbolic skew-products}, i.e., for symbolic skew-products where
the fiber maps are symplectomorphisms.
To do this, let us first recall the method to construct robust
tangencies in symbolic skew-products.

\subsection{Tangencies in symbolic skew-products}
Following~\cite{BR17} we will introduce the notion of tangencies in symbolic
skew-products and afterwards give a criterion for their
construction.

\subsubsection{The set of smooth symbolic skew-products} Since we will need to work with differentiable fiber maps, it will
be useful to extend the set of symbolic skew-products to this
setting.

\begin{defi}
For an integer $r\geq 1$,
$\mathcal{S}^{r}(M)\equiv\mathcal{S}^{r+\alpha}_{\mathscr{A},\nu}(M)$
denotes
 the set of  skew-products $\Phi=\tau\ltimes\phi_\xi$ on $\mathcal{M}=\Sigma\times M$
 such that there are  $\gamma\equiv\gamma(\Phi)>0$,
 $\hat\gamma\equiv\hat\gamma(\Phi)>0$ and
 $C_r\equiv C_r(\Phi)\geq 0$  satisfying
\begin{enumerate}
\item[-] $d_{C^r}(\phi^{\pm 1}_\xi,\phi^{\pm 1}_{\zeta}) \leq
C_r \, d_{\Sigma}(\xi,\zeta)^\alpha$ \ \ for all $\xi,\zeta \in
\Sigma$ with $\xi_0=\zeta_0$, and
\item[-] $\phi_\xi$ are $C^r$-diffeomorphisms of $M$ with
$D^r\phi^{\pm 1}_\xi$ Lipschitz (with uniform Lipschitz constant)
$$
\gamma< m(D\phi_\xi(x)) < \|D\phi_\xi(x)\| <\hat{\gamma}^{-1} \ \
\text{for all $(\xi,x) \in \Sigma\times X$}.
$$
\end{enumerate}
In addition, $
   \mathcal{PHS}^{r}(M)\equiv\mathcal{PHS}^{r+\alpha}_{\mathscr{A},\nu}(M)
   \eqdef \mathcal{PHS}(M)\cap \mathcal{S}^{r}(M)$ for $r\geq 1$.
Finally, a partially hyperbolic skew-product is said to be
\emph{fiber bunched} if $\nu^\alpha<\gamma\hat{\gamma}$.
\end{defi}

We endow $\mathcal{S}^r(M)$  with the metric
$$
 d_{\mathcal{S}^r}(\Phi,\Psi) \eqdef
  d_r(\Phi,\Psi) + \mathrm{Lip}_r(\Phi,\Psi)+
  \mathrm{Hol}_r(\Phi,\Psi)
$$
where
\begin{gather*}
\mathrm{Lip}_r(\Phi,\Psi)\eqdef \max_{\xi \in \Sigma}
  \big|\mathrm{Lip}(D^r\phi_\xi^{\pm 1})-\mathrm{Lip}(D^r\psi_\xi^{\pm
  1})\big| \\
   d_r(\Phi, \Psi)\eqdef \max_{\xi \in \Sigma} \,
 d_{C^r}(\phi^{\pm 1}_\xi,\psi^{\pm 1}_\xi)
\quad  \text{and} \quad
\mathrm{Hol}_r(\Phi,\Psi)\eqdef \big|C_r(\Phi)-C_r(\Psi)\big|.
\end{gather*}
Hence $\mathcal{PHS}^r(M)$ is an open set of $\mathcal{S}^{r}(M)$
and
$ \mathcal{S}^{r+1}(M) \subset \mathcal{S}^{r}(M) \subset
\mathcal{S}(M)$ for
    $r\geq  0$.

\subsubsection{Tangencies} To define the
notion of a tangency for symbolic skew-products we first need to
introduce the notion of a tangent
direction. 


\begin{defi}[tangent direction] Let
$\Phi=\tau \ltimes \phi_\xi\in \mathcal{S}^{r}(M)$ be a symbolic
skew-product with a pair of transitive hyperbolic sets $\Gamma^1$
and $\Gamma^2$  and suppose that $(\xi,x)\in W^{u}(\Gamma^1) \cap
W^{s}(\Gamma^2)$.  A unitary vector $v\in T_x M$ is called a
\emph{tangent direction} at $(\xi,x)$ if there are $C>0$ and
$0<\lambda<1$~such~that
$$
\|D\phi^n_\xi(x)v\| \leq C \lambda^{|n|} \quad \text{for all
$n\in\mathbb{Z}$.}
$$
The maximum number of independent tangent directions at $(\xi,x)$
is denoted by $d_T\equiv d_T(\xi,x)$.
\end{defi}

Now we are ready to give the definition of a tangency.

\begin{defi}[tangency]
We say that $\Phi\in \mathcal{S}^r(M)$ has a \emph{(bundle)
tangency} of dimension $\ell> 0$ between $W^{u}(\Gamma^1)$ and
$W^s(\Gamma^2)$ if there exists $(\xi,x)\in W^u(\Gamma^1)\cap
W^s(\Gamma^2)$ such that
$$
      \ell=d_T(\xi,x) \quad\text{and} \quad   \mathrm{ind}^{cu}(\Gamma^1)+\mathrm{ind}^{cs}(\Gamma^2)-\ell < c.
$$
If $\Gamma^1=\Gamma^2$, the tangency is called \emph{homoclinic},
and otherwise \emph{heteroclinic}. The tangency (of dimension
$\ell$) is said to be \emph{$\mathcal{S}^{r}$-robust} if for any
small enough $\mathcal{S}^{r}$-perturbation $\Psi$ of $\Phi$ has a
tangency (of dimension $\ell$) between the unstable set
$W^{u}(\Gamma^1_\Psi)$ and the stable set $W^s(\Gamma^2_\Psi)$.

The \emph{codimension} of the tangency is defined as
$c_T=c-[\mathrm{ind}^{cu}(\Gamma^1)+\mathrm{ind}^{cs}(\Gamma^2)-\ell]$.
\end{defi}

For the rest of this section, we will work in local
coordinates and thus may assume that $M=\mathbb{R}^c$ with  $c\geq 2$.

\subsubsection{Cone fields in symbolic skew-products}
Consider an integer $1\leq \ell\leq c$. An $\ell$-dimensional
vector subspace of $\mathbb{R}^c$ is called a $\ell$-plane. The
\emph{Grassmannian manifold} $G(\ell,c)$ is defined as the set of
$\ell$-planes in $\mathbb{R}^c$. A standard $\ell$-cone in
$\mathbb{R}^c$ is a set of the form
$$
\mathcal{C}=\{(v,w)\in \mathbb{R}^c: v\in\mathbb{R}^\ell  \
\text{and} \ \|w\| \leq \rho \|v\| \ \ \text{for some} \ \rho>0\}.
$$
More generally, a \emph{$\ell$-cone} is the image of a standard
$\ell$-cone under an invertible linear map. In fact, any
$\ell$-cone $\mathcal{C}$ in $\mathbb{R}^c$ induces an open set in
$G(\ell,c)$, which we will continue denoting by $\mathcal{C}$.

\begin{defi}[stable and unstable cones]
 Let $\Phi=\tau\ltimes\phi_\xi \in
\mathcal{S}^0(\mathbb{R}^c)$ and consider an open set
$\mathcal{B}$ of $\mathcal{M}=\Sigma \times \mathbb{R}^c$. An
$\ell$-cone $\mathcal{C}^{uu}$ in $\mathbb{R}^c$ is said to be
\emph{unstable} for $\Phi$ on $\mathcal{B}$ if there is
$0<\lambda<1$ such that
$$
D\phi_\xi(x)\mathcal{C}^{uu} \subset
\mathrm{int}(\mathcal{C}^{uu})
\quad \text{and} \quad 
\text{$\|D\phi_\xi(x)v\|\geq \lambda^{-1} \|v\|$
 \ \ for all $v\in
\mathcal{C}^{uu}$, $(\xi,x)\in \mathcal{B}\cap
\Phi^{-1}(\mathcal{B})$.}
$$
Similarly we define the \emph{stable $\ell$-cone}
$\mathcal{C}^{ss}$ for $\Phi$ on  $\mathcal{B}$.
\end{defi}

\subsubsection{Tangencies in one-step maps}

Let $\Phi=\tau\ltimes(\phi_1,\dots,\phi_d) \in
\mathcal{PHS}^1(\mathbb{R}^c)$ be a fiber bunched one-step map.
That is,
$$
   \nu^\alpha < \gamma<
   m(D\phi_i)<\|D\phi_i\|<\hat\gamma^{-1}<\nu^{-\alpha} \quad \text{with}
   \quad \gamma\hat\gamma > \nu^\alpha \quad \text{for all $i=1,\dots,d$. }
$$
We fix $0<\ell<c$ and denote
$$
\hat\phi_i(x,E)=(\phi_i(x),D\phi_i(x)E), \ \ \ (x,E)\in
\hat{M}\eqdef \mathbb{R}^c\times G(\ell,c), \quad \text{for
$i=1,\dots,d$}.
$$

In order to provide a criterion to get robust tangencies we need
the following definitions:

\begin{defi}[blending region with tangency]
Let $B$, $D$ be bounded open sets of $\mathbb{R}^c$.
 We say that a $cs$-blending region, $B$, with respect
to $\{\phi_1,\dots,\phi_d\}$ on $D$ has a \emph{tangency} of
dimension $\ell$  if there exist bounded open sets $\hat{B}$ and
$\hat{D}$ in $\hat{M}$ such that
$$
\text{$\hat{B}$ is a $cs$-blending region with respect to
$\{\hat\phi_1,\dots,\hat\phi_d\}$ on $\hat D$ so that $\hat B
\subset B \times \mathcal{C}^{uu}$}.
$$
Here $\mathcal{C}^{uu}$ is an unstable $\ell$-cone  for
$\Phi=\tau\ltimes(\phi_1,\dots,\phi_d)$ on $\Sigma\times B$.  We
say that the set $\hat{B}$ is a \emph{$\ell$-tangency} (of the
blending region $B$). Similarly, we define a $cu$-blending region
with a tangency of dimension $\ell$.
\end{defi}

Let $A_1$ and $A_2$ be two subsets of $\hat{M}$.

\begin{defi}[transition]
\label{def:transition}  We say that  the semigroup $\langle
\hat{\phi}_1,\dots,\hat{\phi}_d\rangle^+$ has a \emph{transition
from $A_1$ to $A_2$} if there exist a point $x\in A_1$ and a map
$T\in \langle \hat\phi_1,\dots,\hat\phi_d\rangle^+$ such that
$T(x)\in A_2$.
\end{defi}

Using these terminologies, the following result from
\cite[Corollary~5.15]{BR17} provides a criterion to construct
robust tangencies:

\begin{thm}
\label{cor:tangencias-IFS} Let
$\Phi=\tau\ltimes(\phi_1,\dots,\phi_d) \in
\mathcal{PHS}^1(\mathbb{R}^c)$ be fiber bunched one-step map with
small enough Lipschitz constant of $D\phi^{\pm 1}_i$. Consider
integers
$$
0<i_1,i_2<c \quad \text{and} \quad \max\{0,i_2-i_1\}<\ell\leq \min
\{c-i_1,i_2\}.
$$
Suppose that there are bounded open sets $D_1, D_2, B_1, B_2$ of
$\mathbb{R}^c$ such that with respect to
$\{\phi_1,\dots,\phi_d\}$,
\begin{enumerate}[leftmargin=0.6cm]
\item[-] $B_1$ is a $cs$-blending region on
$D_1$,  with $cs$-index $i_1$ and a $\ell$-tangency $\hat{B}_1$;
\item[-] $B_2$ is a $cu$-blending region
on $D_2$, with $cs$-index $i_2$ and a $\ell$-tangency $\hat{B}_2$;
\item[-]$\langle\hat\phi_1,\dots,\hat\phi_d \rangle^+$ has a transition from
$\hat{B}_1$ to $\hat{B}_2$.
\end{enumerate}
Then $\Phi$ has a $\mathcal{S}^1$-robust tangency of dimension
$\ell$ between $W^u(\Gamma_1)$ and $W^s(\Gamma_2)$ where
$\Gamma_1$ and $\Gamma_2$ are the maximal invariant sets in
$\Sigma\times \overline{D_1}$ and $\Sigma\times \overline{D_2}$.
\end{thm}

\subsection{Symplectic perturbations}
We would need the perturbative tools stated below. The following
two remarks deal with local perturbations and are done in local
Darboux charts, which are coordinates in which the symplectic form
is written in the canonical way.


\begin{rem}[Pasting lemma]
\label{rem-pasting} A symplectic pasting lemma~\cite[Lemma
3.9]{AM07} states that given a $C^r$-symplectic map with a
periodic point, one can locally around the periodic point, glue
any other $C^r$-close symplectic map. Thus for example, one can
modify the identity map to obtain locally any linear symplectic
matrix.

\end{rem}

The fact described in the next remark (see also~\cite[Proof of
Theorem 3.16: pertubations]{NP12}), says that one can locally
perturb any symplectomorphism by translations.

\begin{rem}[Perturbations by translations]
\label{rm:bump-function} Let $z$ be a point in $M$. Using the
Theorem of Darboux we select local coordinates in a neighborhood
$U$ of $p$ of the form $(x_1,\dots,x_n;y_1,\dots,y_n)$ such that
on $U$ the symplectic form is $dx_1 \wedge dy_1 +\dots + dx_n
\wedge dy_n$. Consider a symplectic diffeomorphism $\phi$ of $M$
and let $U_0$ be a small bounded connected open set so that $z
\in U_0$ and the closure of $\phi(U_0)$ and $U_0$ are contained in
$U$. Given a sufficiently small vector $(u, v) = (u_1, \dots, u_n;
v_1,\dots,v_n)$ in Darboux coordinates, let the Hamiltonian
$H_{(u,v)}$ on $M$ be a bump function such that $H_{(u,v)}(z) = 0$
if $z \in M\setminus U$ and if $z\in U_0$, $H_{(u,v)}$ is
expressed in the local coordinates as $H_{(u,v)}(x,y) = vx - u y$.
Observe that the time-one map $T_{(u,v)}$ of the Hamiltonian flow
$H_{(u,v)}$ in $U_0$ is the translation in local coordinates
$(x,y)\mapsto (x,y)+ (u,v)$ and the identity in $\partial{U}$.
Then, we define the symplectic perturbation of $\phi$ as
$\tilde{\phi}=T_{(u,v)}\circ\phi$.
\end{rem}

Let $\omega_N$ and $\omega_M$ be the symplectic two-forms on the
symplectic manifolds $N$ and $M$. Call $\pi_N$ and $\pi_M$ the
projections of $N\times M$ onto $N$ and $M$. Then a symplectic
form on $N\times M$ is given by
$\omega=\pi_N^\ast\omega_N+\pi_M^\ast\omega_M$ ($\pi^\ast$
indicating the pull-back). The following lemma allows for
semi-local perturbations, and will use the fact that we are
dealing with skew-product symplectomorphisms defined on a product
manifold $N\times M$.

\begin{lem}[semi-local perturbations]
\label{semilocal} Consider $\phi$, a $C^r$-symplectomorphism of
$M$, with the additional hypothesis that $\phi$ is
$C^r$-symplectically isotopic to the identity. Let $F$ be a
$C^r$-symplectomorphism of $N$ and $U_1$, $U_2$ be small enough
neighborhoods contained in some Darboux chart such that
$\overline{U_1} \subset U_2$. Then there exists a
$C^r$-symplectomorphism $f$ of $N\times M$ such that
$f|_{U_1\times M}=F\times \phi$ and $f|_{\partial{U}_2\times
M}=F\times \mathrm{id}$
\end{lem}

\begin{proof}
Let $V$ be the Darboux chart of $N$ so that $U_1, U_2 \subset V$.
We can suppose that $V$ is, in fact, an open ball around the
origin in $\mathbb{R}^{n}$ where $n$ is the dimension of $N$. We
take open balls of radius $r_1<r_2<r_3$ such that $U_1\subset
B_{r_1}(0) \subset B_{r_2}(0) \subset U_{2} \subset B_{r_3}(0)
\subset V$.

Let $\phi:M\to M$ be $C^r$-symplectically isotopic to the
identity, and assume that the isotopy is defined in the interval
$[0,r_3]$, writing $\phi(s,y)$, for $s\in[0,r_3]$. Moreover we can
reparametrize so that $\phi(s,\cdot)=\phi(\cdot)$ for
$s\in[0,r_1]$ and $\phi(s,\cdot)=\mathrm{id}$ for $s\in
[r_2,r_3]$. Consider a map $f$ defined by
$$
f(x,y)=(F(x), \phi(\|x\|, y)) \quad (x,y)\in U_2 \times M.
$$
Notice that $f|_{U_1\times M}=F \times \phi$ and $f|_{\partial
U_2\times M}=F\times \mathrm{id}$. This map can be extended
$C^r$-globally to $N\times M$ simply by defining it as $f=
F\times\mathrm{id}$ outside $U_2\times M$.

It is left to see that $f$ preserves the form $\omega$, and it is
enough to show this locally in product Darboux coordinates
$V\times W\subset N\times M$. We can write the form $\omega_N$ in
$V$ as $\sum dx_i\wedge dx_j$ and the form $\omega_M$ in $W$ as
$\sum dy_i\wedge dy_j$, and then $\omega=\sum dx_i\wedge dx_j+\sum
dy_k\wedge dy_\ell$ where $(x,y)=(x_1,\dots,x_{n},y_{n+1},\dots,
y_{n+c})$ are the coordinates in $V\times W$. Then
$$
f^\ast\omega=\sum \det\big(\frac{\partial f_i}{\partial x_j}\big)
\ dx_i\wedge dx_j+ \sum \det\big(\frac{\partial f_k}{\partial
y_\ell}\big) \ dy_k\wedge dy_\ell.
$$
Since $F$ is  symplectic then $\det(\frac{\partial f_i}{\partial
x_j})=\det(\frac{\partial F_i}{\partial x_j})=1$. For the other
hand,
$$
\det\big(\frac{\partial f_k}{\partial
y_\ell}\big)=\det\big(\frac{\partial \phi_k}{\partial
y_\ell}\big), \quad \text{where the index $k>n$.}
$$
Since for any $x$ fixed, the map $\phi(||x||, \cdot)$ preserves
the form $\sum dy_k\wedge dy_\ell$, this means that writing
$D\phi=(\partial_x \phi,
\partial_y \phi)$, the submatrix  $\partial_y\phi$ is symplectic.
Going back this implies that $\det(\frac{\partial f_k}{\partial
y_l})=1$, when $k>n$, and thus $f^\ast\omega=\omega$.
\end{proof}

\subsection{Tangencies and transitivity for symplectic symbolic skew-products}
In this subsection we want to get the following result:

\begin{thm}
\label{thm:tg}
Let $M$ be a symplectic manifold of dimension $c\geq 2$ and consider
a 
integer $0<\ell\leq c/2$. Then, there are $d=d(c)\geq
2$ and an arc of symplectic one-step maps 
$\Phi_\epsilon=\tau\ltimes (\phi_1, \dots, \phi_d)$ in
$\mathcal{PHS}^r(M)$ isotopic to $\Phi_0=\tau\times \mathrm{id}$
so that for any $\varepsilon>0$
\begin{enumerate}
\item[-] $\Phi_\varepsilon$  is $\mathcal{S}^0$-robustly
topologically mixing;
\item[-] $\Phi_\varepsilon$ has a
$\mathcal{S}^1$-robust  tangency of dimension
$\ell$. 
\end{enumerate}
\end{thm}

In order to prove the above theorem, according to
Theorem~\ref{cor:tangencias-IFS} and Corollary~\ref{cor:transitividad}
we need to create arcs of IFSs
generated by symplectomorphisms having blending regions with
tangency, transitions and globalization. Namely, fixing an integer $0<\ell\leq
c/2$, we need to prove the following result:

\begin{prop}
\label{key} There are  $d=d(c)\geq 3$ and arcs of
$C^r$-symplectomorphisms of $M$
$$
\phi_1\equiv\phi_1(\varepsilon), \dots,
\phi_d\equiv\phi_d(\varepsilon), \ \ \varepsilon\geq 0,  \quad
\text{with $\phi_i(0)=\mathrm{id}$}
$$
and having the following properties: for any $\varepsilon>0$,
there are bounded open sets $B_1\equiv B_1(\varepsilon)$ and
$B_2\equiv B_2(\varepsilon)$ of $M$ such that with respect to
$\{\phi_1,\dots,\phi_d\}$,
\begin{enumerate}
 \item $B_1$ is a $double$-blending region with $cs$-index $\ell$ and having a $\ell$-tangency
 $\hat{B}_1$;
 \item $B_2$ is a $cu$-blending region with $cu$-index $\ell$ and having a $\ell$-tangency
 $\hat{B}_2$;
 \item $\langle \hat\phi_1,\dots, \hat\phi_d\rangle^+$ has transition from $\hat{B}_1$ to $\hat{B}_2$ and/or transition from
$\hat{B}_1$ to $\hat{B}_1$;
\item $B_1$ is globalized by $\langle \phi_1,\dots, \phi_d\rangle^+$.
\end{enumerate}
\end{prop}

Now let us show the existence of such symplectic arcs. We will
always work locally inside Darboux charts. To create local
blending regions $B_1$ and $B_2$, one can take a symplectic map
with a local fixed hyperbolic point and compose the map with the
necessary local translations as described in the
Remark~\ref{rm:bump-function} where the directions of the
translations come from Proposition~\ref{prop:creation-blender}.
The hardest work is to obtain a $\ell$-tangency for these local
blending regions. First we need some basic facts from symplectic geometry:

A \emph{symplectic vector subspace} is one on which the symplectic
form still restricts to a non-degenerate two-form. By a
\emph{coisotropic vector subspace} is understood a subspace of
$\mathbb{R}^c$ which contains its symplectic complement.

\begin{lem}
\label{lema:tang} Let $\phi$ be a hyperbolic linear map of
$\mathbb{R}^c$ that preserves the symplectic form, i.e., a
symplectic matrix, having a splitting in eigenspaces at the fixed
point $p=0$ of the form $E^s\oplus E^{cu}\oplus E^{uu}$ where the
unstable direction is $E^u=E^{cu}\oplus E^{uu}$ and $0<\ell=\dim
E^{uu}<c$.  If $\ell$ is odd we ask that $E^{uu}$ is a coisotropic subspace and
a symplectic subspace if $\ell$ is even.

Then, there are
\begin{enumerate}
\item an unstable $\ell$-cone $\mathcal{C}^{uu}$ around $E^{uu}$,
\item neighborhoods $B$ of $p$ in $\mathbb{R}^c$ and $G\subset \mathcal{C}^{uu}$ of $E^{uu}$ in
$G(\ell,c)$, and
\item symplectic ortogonal matrizes $A_1,\dots,A_d$ and
vectors $c_1,\dots,c_d$ in $\mathbb{R}^c$ with $d=d(c)\geq 2$
\end{enumerate}
such that
$$
\overline{\hat B}\subset\bigcup_{i=1}^d \hat\phi_i(\hat B) \quad
\text{where} \quad \hat B = B \times G
$$
and $\hat\phi_i$ are the induced maps on
$\hat{M}=\mathbb{R}^c\times G(\ell,c)$  by the symplectomorphisms
$\phi_i=A_i\cdot \phi+c_i$.  That is, $B$ is a $cs$-blending
region with respect to $\{\phi_1,\dots,\phi_d\}$ of $cs$-index
$\dim E^s$ and a $\ell$-tangency. Moreover, the maps
$\phi_1,\dots,\phi_d$ can be constructed as an arc of
symplectomorphisms homotopic to $\phi$.
 \end{lem}
\begin{proof}
 As in Proposition~\ref{prop:creation-blender}, we can take
 a small neighborhood $B$ of $p$ in $\mathbb{R}^c$
 and choose arcs of translations $f_i=\phi+c_i$ homotopic to $\phi$ with
 $c_i\in\mathbb{R}^c$ in order to have that $B$ is
 a $cs$-blending region with respect to $\{f_1,\dots,f_{k_1}\}$. In
 particular, $\{ f_i(B): i=1,\dots,k_1\}$ is an open cover of
 $\overline{B}$. Notice that there exists an $\epsilon>0$,
 such that the above covering property holds for any
 $\epsilon$-perturbation $\phi_i$ of $f_i$.

 On the other hand, $D\phi(p)$ induces a map
 $A$ on $G(\ell,c)$ given by $A(E)=D\phi(p)E$.
 Since $\phi$ is $C^2$, this map
 has a hyperbolic attracting fixed point $E^{uu}$.
 Then, as in Proposition~\ref{prop:creation-blender},
 there are arcs $T_1\equiv T_1(\varepsilon),\dots,T_{k_2}\equiv T_{k_2}(\varepsilon)$
 of translations on $G(\ell,c)$ homotopic to the identity and
 a neighborhood $G\equiv G(\varepsilon)$ of
 $E^{uu}$ in $G({k,c})$, such that
 $
 \{F_j(G): j=1,\dots,k_2\}$
 is
 an open cover of $\overline{G}$ where $F_j= T_j \circ A$.
 Each translation map $F_j$ in the Grassmannian corresponds to a
map of the form $A_j\cdot D\phi(p)$ in the tangent bundle where
$A_j$ is an orthogonal matrix and $\|A_j- \mathrm{id}\|\to 0$ as
$\varepsilon\to 0$ for all $j=1,\dots,k_2$. Moreover, we can take
$G$ small enough so that $G\subset \mathcal{C}^{uu}$ where
$\mathcal{C}^{uu}$ is a unstable $\ell$-cone around $E^{uu}$ for
$\phi$. Taking
$$\|A_j-\mathrm{id}\|< \epsilon / \|\phi\| \quad \text{and}\quad
\phi_{ij}=A_j \cdot \phi + c_i$$ then $\|\phi_{ij}-f_i\|<\epsilon$
and so
$$
\overline{B} \subset \bigcup_{i=1}^{k_1}
\bigcup_{j=1}^{k_2}\phi_{ij}(B).
$$

The translations $c_i$ preserve the symplectic form, and so we  only need
to show the matrix $A_j$ can be taken symplectic.
Suppose for now that $\ell$ is even and so $E^{uu}$ is a symplectic
subspace.
Observe that the symplectic subspaces form an open
and dense subset in the set of $\ell$-planes.
Then we can take the $\ell$-cone $\mathcal{C}^{uu}$
around $E^{uu}$ in $G(\ell,c)$ as
an open set of symplectic subspaces. Since $A_j$ is arbitrarily close to the
identity, we may assume that $E=A_j E^{uu}$ is  in $\mathcal{C}^{uu}$
and so is also a symplectic subspace. As symplectic
matrices act transitively on symplectic vector
subspaces~\cite[Cap.~V]{Zehnder}, we can take a symplectic
matrix $\tilde{A}_j$ , close to the identity,  such that $\tilde{A}_j E^{uu}=E$.
Without loss of generality we assume $\tilde{A}_j=A_j$.
If $\ell$ is an odd number, we proceed similarly.
Since coisotropic subspaces are an
open an dense set of $G(\ell,c)$ where  the symplectic
matrices also act transitively~\cite[Exer.~2.33]{MS98}, then in the same
manner as before, we can take $A_j$  as a symplectic matrix.

Finally, by construction, the maps $\phi_{ij}$ induce a set of
maps $\hat{\phi}_{ij}$ on $\hat{M}$ so that
$\{\hat{\phi}_{ij}(\hat{B})\}$ is an open cover of the closure of
$\hat{B}=B\times G$. Moreover, since $A_j$ is an orthogonal matrix
then the hyperbolicity of the blending regions still holds. That
is, $B$ is also a $cs$-blending region with respect to
$\{\phi_{ij}: i=1,\dots,k_1, \ j=1,\dots,k_2\}$ with $cs$-index
equal to $\dim E^s$ and as well is a $\ell$-tangency (the
$cs$-blending region $\hat{B}$ on $\hat{M}$). This completes the
proof.
\end{proof}

Notice that, by means of an arbitrarily small perturbation of
the identity map, we can create a map $\phi$ for which $x$ is a
hyperbolic fixed point. Thus the above construction can be done
homotopic to the identity.

\begin{proof}[Proof of Proposition~\ref{key}]
 To make the transition map from  $\hat{B}_1$ to $\hat{B}_2$,
the blending regions can be chosen inside the same Darboux chart
in such a manner so that the transition map is simply a
translation. On the other hand  the transition map from
$\hat{B}_1$ to $\hat{B}_1$, can be chosen as a local linear
symplectic map that has a fixed point inside the double-blending
region and takes the unstable cone into the stable one over
several iterates, or for example maps $E^{uu}$ into $E^{ss}$.
Recall that these subspaces can be chosen to be symplectic or
coisotropic on which the symplectic matrices act
transitively~\cite[Cap.~V]{Zehnder} and~\cite[Exer.~2.33]{MS98}.

To globalize the blending region $B_1$,
observe that the $C^r$-diffeomorphisms constructed
to prove Theorem~\ref{thmD} are locally small translations. Thus
again using Darboux coordinates and a
Hamiltonian bump function as in the Remark~\ref{rm:bump-function}, one can glue
the maps required for the Lemma~\ref{lem:Rn} and
Proposition~\ref{thm:globalizacion-id}. The Milnor atlas required for the
globalization has to be taken subordinate to the Darboux atlas.
\end{proof}

The proof of Theorem~\ref{thm:tg} is now complete.

\subsection{Symplectic realization: proof of Theorem~\ref{thmD}}
Recall that our goal is to construct arcs of
$C^r$-symplectomorphisms $f_{\varepsilon}$ of $N \times M$
isotopic to $f_{0} = F \times \mathrm{id}$ having robust
homoclinic (or equidimensional) tangencies of codimension
$\ell\leq\dim(M)/2$. From Theorem~\ref{thm:tg}, suppose
$\Phi_{\varepsilon}=\tau\ltimes(\phi_1,\ldots,\phi_d)\in
\mathcal{PHS}^{1}(M)$ is an arc of robustly transitive symplectic
symbolic skew-products with a robust homoclinic tangency of
dimension $\ell$. It is desirable to obtain a one-parameter family
of
diffeomorphisms $f_\varepsilon$ satisfying 
$ f_{\varepsilon}|_{R_i\times M}=F \times \phi_i$ 
for $i=1,\dots,d$  where $R_i$ is the Markov partition of the set
$\Lambda$. The Markov partition must be small enough so that each
rectangle $R_i$ is inside the Darboux chart $U_i$. If the initial
set $\Lambda$ is conjugated to a shift of $d$ symbols, the pieces
of the new partition would correspond to cylinders $\mathsf{C}_i$
in the shift. A cylinder consists of all sequences with a given
block of indices fixed around the zeroth coordinate, the length of
the cylinder being the size of the given block. With each $U_i$ we
associate the map $\phi_{k_i}$  where the index $k_i$ corresponds
to the zeroth index of the cylinder $\mathsf{C}_i$.

Since $U_i$ are disjoint, we may apply Lemma~\ref{semilocal}
inductively with respect to the maps $\phi_{k_i}$, and thus obtain
a symplectomorphism $f_\varepsilon$ in $N\times M$ satisfying
\begin{align*}
f_{\varepsilon} |_{R_i\times M}&=F \times \phi_{k_i} \quad \text{
$f_{\varepsilon}=F \times \mathrm{id}$  outside of $\cup
(U_i\times M)$}.
\end{align*}
Then by construction the dynamics over $\Lambda\times M$ will be
conjugated to $\Phi_\varepsilon$. Using this conjugation one can obtain
robust tangencies and transitive sets for the map $f_\varepsilon$, exactly as was
proven in~\cite[Sec.~6.1]{BR17}. This concludes the proof of
Theorem~\ref{thmD}.

